\documentclass{article}
\usepackage{amssymb,amsmath}

\input xy
\xyoption{all}

\begin{document}
\title{Cohomology of the de Rham complex twisted by the oscillatory representation}
\author{Svatopluk Kr\'ysl \footnote{{\it E-mail address:} Svatopluk.Krysl@mff.cuni.cz, {\it Tel./Fax:} $+$ 420  222 323 221/$+$ 420 222 323 394} \\{\it \small 
Charles University in Prague, Sokolovsk\'a 83, Praha, 186 75,  Czech Republic}
\thanks{I thank to Mike Eastwood for introducing me into the world of representation theory of Lie groups and its applications in differential geometry, and for his numerous 
inspiring lectures.}
}
\maketitle \noindent
\centerline{\large\bf Abstract}

We introduce a Hilbert $A$-module structure on the higher oscillatory module, where $A$ denotes the  $C^*$-algebra  of bounded endomorphisms
of  the basic oscillatory module. We also define the notion of an exterior covariant derivative in an $A$-Hilbert bundle and
use it for a construction of an $A$-elliptic complex of differential operators for certain symplectic manifolds equipped with a symplectic connection.

\section{Introduction}
 
In mathematical physics, one often uses complexes of exterior forms twisted by finite rank vector bundles.
For example, the non-quantum gauge theory of the electromagnetic interaction can be seen as the mathematics of 
connections on line bundles over four dimensional   manifolds. The connections are models of the abelian $U(1)$-gauge fields.
Notice also that the Maxwell equations for vacuo can be formulated by taking the line bundle connection to be trivial. 
Thus, one investigates the de Rham complex over Lorentzian manifolds.
In Riemannian geometry, a similar sequence, namely
$(d_i^{\nabla^{g}},\Gamma(\bigwedge^{i}T^*M\otimes TM))_{i \in \mathbb{N}_0}$ can be considered. 
Here, $d_i^{\nabla^{g}}$ denotes the exterior covariant derivative induced
by the Riemannian connection $\nabla^{g}$ of the Riemannian manifold $(M^n,g).$ 
This sequence  might be used to treat, e.g., the $SO(n)$-invariant parts of the curvature operator and the Bochner-Weitzenb\"ock equations systematically.
However, non-flat manifolds, for which the sequence fails to be a complex, are more in the focus of Riemannian geometers.
Roughly speaking, the Riemannian counterpart of the Maxwell equations does not seem to be so useful and widely used. 
  
 In symplectic geometry, the role of an ``adapted`` connection is not as obvious as in the Riemannian geometry, foremost because of  a theorem of Darboux due to which all
symplectic manifolds are equivalent at the infinitesimal level.
Therefore it is not surprising that for any symplectic manifold $(M,\omega),$ the torsion-free symplectic  connections
form an infinite dimensional affine space. See, e.g., Tondeur \cite{T} or Gelfand et al. \cite{GRS}. The role of  symplectic  
 connections seems to be rather in  
Quantization of classical Mechanics. See Fedosov \cite{Fedosov} and notice that the  higher oscillatory module studied  here
is related to the Weyl algebra structure used in his quantization procedure.

In Fomenko, Mischenko \cite{FM}, the concept of  $A$-Hilbert bundles and differential operators for them were introduced for a unital $C^*$-algebra $A.$ These bundles have the so called Hilbert $A$-modules  as fibers.  
See   Solovyov, Troitsky \cite{ST}   for these notions.
The authors of \cite{FM} investigate  finitely generated projective $A$-Hilbert bundles over compact manifolds and 
$A$-elliptic operators acting between their smooth sections. They prove that such operators have the so called $A$-Fredholm property. In particular, their
 kernels are finitely generated projective Hilbert $A$-modules.
In Kr\'ysl \cite{Krysl}, the results of \cite{FM} were 
transferred to $A$-elliptic complexes and conclusions for the cohomology groups of these complexes were made.  

In this paper, we introduce a sequence  of infinite rank vector bundles over a symplectic manifold and differential operators acting
between sections of these bundles. The symplectic manifold $(M,\omega)$ is supposed to admit a {\it metaplectic structure},
 a symplectic analogue of the Riemannian spin structure, and a symplectic connection $\nabla.$
The principal group  of the metaplectic bundle  is the metaplectic group $Mp(2n,\mathbb{R}).$
At least since the articles of Shale \cite{Shale} and Weil \cite{Weil},
a faithful unitary representation of the metaplectic group on the space $H=L^2(\mathbb{R}^n)$ is known. This is 
the so called Segal-Shale-Weil, metaplectic or symplectic spinor representation.
In this text, however, we call this representation
the  {\it basic oscillator module}  stressing the fact that it is used as a state space of the 
quantum  harmonic oscillator.   
Associating the basic oscillatory module to the metaplectic structure, we get the so called basic oscillator or symplectic
spinor bundle denoted by $\mathcal{H}.$ See Kostant \cite{Kostant} and Habermann, Habermann \cite{KH}.
The sequence of bundles mentioned at the beginning of this paragraph is formed by the tensor products of the bundle $\bigwedge^{\bullet} T^kM$ of exterior differential $k$-forms on $M$ and of $\mathcal{H}.$
The lift of the symplectic connection $\nabla$ to the basic oscillatory bundle $\mathcal{H}$ induces  the exterior covariant derivatives $d_i^{\nabla^H}: \Gamma(\bigwedge^{i}T^*M \otimes \mathcal{H}) \to 
\Gamma(\bigwedge^{i+1}T^*M \otimes \mathcal{H}).$
If the curvature of $\nabla$ is zero, the sequence $d_H^{\bullet}=(d_i^{\nabla^H},\Gamma(\bigwedge^{i}T^*M \otimes \mathcal{H}))_{i \in \mathbb{N}_0}$ forms a complex. 
This the de Rham complex twisted by the oscillatory module. We prove that this complex is $A$-elliptic and use a result from \cite{Krysl} to get an information on the cohomology groups
of this complex when $M$ is compact. As far as we know, this is the first $A$-elliptic complex in infinite rank vector bundles described in the literature.

In the second section, we recall the notion of a Hilbert $A$-module and introduce  the higher oscillatory module as a module over the metaplectic group as well as over the 
unital  $C^*$-algebra $A$ of bounded endomorphism of $L^2(\mathbb{R}^n).$
We prove   that the oscillatory module is a finitely generated projective Hilbert $A$-module with respect to 
a natural $A$-product (Lemma 3). In section 3, we recall the definition  of an $A$-Hilbert bundle and introduce the notion of
the exterior covariant derivative in an $A$-Hilbert bundle and compute its symbol (Theorem 5).
In  section 3.1., we give a definition of an $A$-elliptic complex, construct the de Rham complex twisted by the basic oscillatory module and state 
a theorem on the properties of its cohomology groups of this complex (Theorem 6).

\section{Higher oscillatory modules}

Let $A$ be a  $C^*$-algebra.
Let us recall a definition of an $A$-Hilbert module briefly. 
For general $C^*$-algebras, this notion was firstly considered by Paschke in \cite{Paschke}. A pre-Hilbert $A$-module  is  any left $A$-module $U$ equipped with a map
$(,):U \times U \to A$ satisfying for each $r\in \mathbb{C},$ $u,v,w \in U$ and $a\in A$
\begin{itemize}
\item[1)] $(ru+v,w)=\overline{r}(u,w)+(v,w)$ 
\item[2)] $(a.u,v)=a^*(u,v)$
\item[3)] $(u,v)=(v,u)^*$ 
\item[4)] $(u,u) \geq 0$ and $(v,v)=0$ implies $v=0.$ 
\end{itemize}

The relation $a\geq b$ holds for $a,b \in A$ if and only if $a-b$ is hermitian and its spectrum lies in $\mathbb{R}^+_0$.
Notice that from 2) and 3), we get $(u,a.v)=(u,v)a$ for any $a\in A$ and $u,v \in U.$
A pre-Hilbert $A$-module is called {\it Hilbert $A$-module} if it is a complete space with respect to the norm $|\,|_U:U \to \mathbb{R}$
defined by $|u|= \sqrt{|(u,u)|_A}.$  In particular, a Hilbert $A$-module  is a Banach space.  If $U$ is a pre-Hilbert $A$-module, we speak
of $(,)$ as of an $A$-product. When $U$ is a Hilbert $A$-module, we call the $A$-valued map $(,)$ a Hilbert $A$-product.

 In the category of pre-Hilbert $A$-modules, the set of morphisms is formed by continuous $A$-equivariant maps between the objects.
Continuity is meant with respect to the (possibly non-complete)
norms. Declaring  the category of Hilbert $A$-modules to be a full subcategory of the category of pre-Hilbert $A$-modules, defines the set of 
morphisms in this category.
 Let us  notice that adjoints  are considered with respect to the $A$-products, i.e., for a morphism $B: U \to V$ of pre-Hilbert $A$-modules,
its adjoint, denoted by $B^*,$ is the map $B^*: V \to U$ such that for each $u\in U$ and $v\in V,$ the equation $(Bu,v)=(u,B^*v)$ holds. 
It is known  that unlike for Hilbert spaces, morphisms of Hilbert $A$-modules do not have adjoints in general. 
When we write a direct sum of Hilbert $A$-modules, we mean that the summands are mutually orthogonal  with respect to $(,).$ In general,  orthocomplements do
 not have the ``exhaustion`` property, i.e., there exist a Hilbert $A$-module $U$ and a (closed) Hilbert $A$-submodule 
$V$ of $U$ such that $U \neq V \oplus V^{\perp}.$ 
(See, e.g., Lance \cite{Lance} for an example.) Fortunately, we have the following

{\bf Theorem 1:} Let $U,V$ be Hilbert $A$-modules and $B: U \to V$ be a Hilbert $A$-module morphism. If the adjoint of $B$ exists and $\mbox{Im }B$ is closed, then
$U  = \mbox{Ker }B \oplus \mbox{Im }B^*.$

{\it Proof.} See  Lance \cite{Lance} (Theorem 3.2) for a proof. $\Box$

Now, we focus our attention to the higher oscillatory module.
Let $(V,\omega)$ be a real symplectic vector space of dimension $2n$ 
and $g:V \times V \to \mathbb{R}$ be a scalar product on $V$. For any $\xi \in V^*,$ we denote by $\xi^g$ the vector in $V$ defined by $\xi(v)=g(\xi^g,v)$ for any $v\in V.$
We denote the appropriate 
extension of $g$ to $\bigwedge^{\bullet}V^*$  by $g$ as well. Let $\tilde{G}$ be a realization of the metaplectic group associated to the symplectic space and 
let $\lambda$ be the  covering homomorphism of the symplectic group of $(V,\omega)$
by the metaplectic group $\tilde{G}.$

Further, we  denote the exterior multiplication of exterior $k$-forms by a $1$-form $\xi$ by $\mbox{ext}^{\xi}_k$ and   recall the following lemma (of \'E. Cartan)
 usually proved by induction
on the dimension of $V$. We use this lemma when we will investigate the $A$-ellipticity of the de Rham complex twisted by the basic oscillatory module.

{\bf Lemma 2:} For any $\xi \in V^* \setminus \{0\},$ the complex $\mbox{ext}^{\bullet}_{\xi} = (\mbox{ext}_i^{\xi}, \bigwedge^i V^*)_{i=-1}^{2n+1}$ is exact.

{\it Proof.} Consider $\bigwedge^i V^*$ as the wedge power of the defining representation of $SO(V,g)$ and suppose $\xi \in V^*.$
Because $\xi\wedge \xi \wedge \alpha = 0$ for any $\alpha \in \bigwedge^k V^*,$ $\mbox{ext}^{\bullet}_{\xi}$ is a complex.
Now, suppose $\xi \wedge \alpha = 0$ for a $k$-form  $\alpha.$ Making contraction of this equation by $\xi^g,$ we get 
$(\iota_{\xi^g}\xi) \alpha - \xi \wedge \iota_{\xi^g}\alpha = 0.$ From that
$\alpha = -(g(\xi,\xi))^{-1} \xi \wedge \iota_{\xi^g}\alpha$ provided
$g(\xi,\xi)\neq 0,$ which holds if and only if $\xi \neq 0.$ Thus $\alpha \in \mbox{Im} (\mbox{ext}_i^{\xi}).$
$\Box$

Let $L$ be a Lagrangian subspace of $(V,\omega).$ 
The Segal-Shale-Weil representation (SSW-representation) of the metaplectic group $\tilde{G}$  is a faithful unitary representation
of this group on the complex Hilbert space $H=L^2(L,g_{|L \times L}).$ By writing $L^2(L, g_{|L\times L}),$ we stress the fact that the metric structure on $L$ is fixed.
Let us denote the SSW-representation by $\rho_0:\tilde{G} \to \mbox{Aut}(H)$ and the scalar product on $H$ by $(,)_H,$ i.e.,
$$(k,l)_H=\int_{x\in L} \overline{k(x)} l(x) dx \, \mbox{ for each } k, l \in H.$$
See Shale \cite{Shale}, Weil \cite{Weil} and Kashiwara, Vergne \cite{KV} for more information on $\rho_0.$

Let us consider the tensor product $\rho$ of  the wedge powers of the dual of the representation $\lambda: \tilde{G} \to \mbox{Aut}(V)$ and the SSW-representation $\rho_0$, i.e.,
we consider a representation $\rho: \tilde{G} \to \mbox{Aut}(C^{\bullet})$ of the metaplectic group on the space $C^{\bullet}=\bigwedge^{\bullet}V^*\otimes H.$ 
Because $\bigwedge^{\bullet}V^*$ is finite dimensional, 
we may say that we consider the canonical Hilbert space topology on the space $C^{\bullet}=\bigwedge^{\bullet} V^* \otimes H.$ Note that $\rho$ is not unitary unless $V=\mathbb{R}.$
The $\tilde{G}$-module $C^{\bullet}$ is the {\it higher oscillatory module}. We call $C^0=H$ the {\it basic oscillatory module}.  
Notice that it is known that the basic oscillatory module splits into two irreducible representations of
$\tilde{G}$ namely, into to the space of even and odd square integrable functions on $L.$ For the
decomposition of $C^{\bullet}$ into irreducible $\tilde{G}$-submodules see, e.g., Kr\'ysl \cite{KryslJOLT2}.

Now, we would like to investigate $C^{\bullet}$ from an "analytical" point of view. Let $A=\mbox{End}(H)$ be
 the unital  $C^*$-algebra of continuous endomorphisms of $H.$ The star operation $*:A \to A$ in $A$
is given by the adjoint of endomorphisms, i.e., $*a = a^*$ for any $a\in A.$ As the norm in $A,$ we take the classical supremum norm, i.e., for any $a\in A,$ we set 
$|a|_A = \mbox{sup}_{|k|_H\leq 1} |a(k)|_H,$ where $| \,|_H$ denotes the norm on $H$ derived from the scalar product $(,)_H.$
Let us remark, that we consider everywhere defined operators only. In particular, the star $*$ is a well defined (and continuous) anti-involutive map in the Banach algebra $A.$

 The space $C^{\bullet}$ introduced above is not only a $\tilde{G}$-module, but it is also an $A$-module upon the action
$$a.(\alpha \otimes k)=\alpha \otimes a(k), \mbox{  } \alpha\otimes k \in C^{\bullet} \, \mbox{and } a\in A.$$ 
For any $k\in H,$ let $k^*: H \to \mathbb{C}$ denote the $(,)_H$-dual to $k$, i.e., $k^*(l)=(k,l)_H.$
Now, let us introduce a product $(,)$  with values in $A$ on the higher oscillatory module 
$C^{\bullet}.$  For any $\alpha \otimes k, \beta \otimes l \in C^{\bullet}$ we set $$(\alpha \otimes k, \beta \otimes l) = g(\alpha, \beta) k\otimes l^* \in A$$ where by
$k\otimes l^*$ we mean the element of $A$ given by $(k\otimes l^*)(m)=l^*(m) k \in H$ for all $m\in H.$ The product is extended to non-homogeneous elements  linearly.
Let us make the following   observation which we use later.
 For $k,l \in H,$ we have $$(k \otimes l^*)^* = l \otimes k^*.$$ Indeed, for any $m,n \in H,$	
we have $((k \otimes l^*)^*m,n)_H=(m, (k\otimes l^*)n)_H=(m, k(l,n)_H)_H= (l,n)_H(m,k)_H = (\overline{(m,k)_H}l,n)_H =
(l(k,m)_H,n)_H = ((l\otimes k^*)m,n)_H.$

{\bf Lemma 3:} The space $C^{\bullet}$ together with $(,)$ is a finitely generated projective Hilbert $A$-module.

{\it Proof.} 
Let $e_0$ be a unit length vector in $H$ and $v$ an arbitrary element of $H.$ The map
$b=  v\otimes e_0^*$ has the property $b(e_0) = v$ and $|b|_A\leq |v|_H.$ 
Let $\{U_i\}_{i=1}^{2^{2n}}$ be an orthonormal basis of $\bigwedge^{\bullet} V^{*}.$ Then obviously, $(U_i \otimes e_0)_{i=1}^{4^{n}}$ is a set of generators of $C^{\bullet}.$
Thus, $C^{\bullet}$ is finitely generated over $A.$

Now, we prove that $C^{\bullet}$ is a Hilbert $A$-module.
\begin{itemize}
\item[1)] $A$-linearity of $(,).$ For any $a\in A,$ $\alpha\otimes k, \beta \otimes l \in C^{\bullet}$ and $m\in H,$ we have 
\begin{eqnarray*}
(a.(\alpha \otimes k), \beta \otimes l)(m) & = & (\alpha \otimes a(k), \beta \otimes l)(m)\\
&=&g(\alpha, \beta) (a(k) \otimes l^*)(m)\\
&=&g(\alpha,\beta)a(k)(l, m)_H
\end{eqnarray*}
On the other hand, we have
\begin{eqnarray*}
a(\alpha \otimes k, \beta \otimes l)(m) & = & g(\alpha,\beta) a(k \otimes l^*)(m)\\
&=&a(k) g(\alpha, \beta) (l,m)_H
\end{eqnarray*}

\item[2)] Self-adjointness. For any $\alpha\otimes k, \beta \otimes l \in C^{\bullet},$ we have
\begin{eqnarray*}
(\alpha \otimes k, \beta \otimes l)^*&=& g(\alpha,\beta) (k \otimes l^*)^*=g(\alpha,\beta) (l\otimes k^*)\\
&=&g(\beta,\alpha)(l\otimes k^*)= (\beta \otimes l, \alpha \otimes k)
\end{eqnarray*}

\item[3)] Positive definiteness.
 Let $c = \sum_{i=1}^{4^{n}} U_i \otimes k_i^*$ for $k_i \in H,$ $i=1,\ldots, 4^{n}.$
Then $(c,c)=\sum_{i,j=1}^{4^{n}} g(U_i,U_j)(k_i\otimes k_j^*)=\sum_{i,j}^{4^{n}} \delta_{ij}(k_i \otimes k_j^*)=
\sum_{i=1}^{4^{n}}(k_i \otimes k_i^*).$ The spectrum of each of the summands consists of the non-negative number
$(k_i,k_i)_H.$ Thus, $k_i \otimes k_i^{*} \geq 0$. Because the non-negative elements in a  $C^*$-algebra form a cone, $c$ satisfies $c \geq 0.$
Suppose $(c,c)_H = 0$ and that the summand in $\sum_{i=1}^{4^{n}}(k_i \otimes k_i^*)$ with index $i_0$ is non-zero.
Writing $-(k_{i_0}\otimes k_{i_0}^*) = \sum_{i \in \{1,\ldots, 4^{n}\}\setminus\{i_0\}} k_{i} \otimes k_{i}^*$ gives a contradiction.
 
\item[4)] Completeness is obvious because the normed space $C^{\bullet}$ is a finite sum of copies of the Hilbert spaces $H$ and the  norm
is derived from the Hilbert scalar products on these copies.
\end{itemize}

Since as we already know $C^{\bullet}$ is a finitely generated Hilbert $A$-module over a unital  $C^*$-algebra, it is projective (see Frank, Larsen \cite{FL},  Theorem 5.9).
$\Box$

\section{Covariant derivatives  and the twisted de Rham complex}

  Let $M$ be a  manifold and $p: \mathcal{E} \to M$ be an $A$-Hilbert bundle, where
$A$ is a fixed unital  $C^*$-algebra. This means in particular, that $p$ is a smooth Banach bundle the fibers of
which are isomorphic to a fixed Hilbert $A$-module $U.$
As it is standard, we denote the space of smooth sections of $\mathcal{E}$ by $\Gamma(\mathcal{E}).$ For any $m\in M,$ the fiber $p^{-1}(m)$ is denoted by $\mathcal{E}_m.$
The morphisms between $A$-Hilbert bundles $p_i:\mathcal{E}_i \to M,$ $i=1,2,$ are supposed to be smooth bundle maps $S: \mathcal{E}_1 \to \mathcal{E}_2,$ i.e.,
$p_1 = p_2 \circ S$ such that for each point $m \in M,$ $S_{|(\mathcal{E}_1)_m}:(\mathcal{E}_1)_m \to (\mathcal{E}_2)_m$ is a morphism of $A$-Hilbert modules.
See, Solovyov, Troitsky \cite{ST} for more information on $A$-Hilbert bundles. 

Let us denote the trivial $A$-Hilbert bundle $M\times A \to M$ by $\mathcal{A}$ and
introduce the product $(,)_{\mathcal{A}}: \Gamma(\mathcal{E}) \times \Gamma(\mathcal{E}) \to \Gamma(\mathcal{A})=\mathcal{C}^{\infty}(M,A)$
on $\Gamma(\mathcal{E})$ by the formula
$$(s,t)_{\mathcal{A}}(m)=(s(m),t(m))_m \in A,$$ where $s,t\in \Gamma(\mathcal{E}),$ $(,)_m$ is the Hilbert $A$-product in $\mathcal{E}_m$ and $m\in M.$ 

Now, let us choose a Riemannian metric $g$ on $M$ and denote by   $|\mbox{vol}_g|$ a choice of the volume element associated to $g.$
 From now on, we suppose that $M$ is compact.
The space $\Gamma(\mathcal{E})$ of smooth sections of $p:\mathcal{E} \to M$ carries a pre-Hilbert $A$-module structure.
The action of $A$ on $\Gamma(\mathcal{E})$ is defined by $(a.s)(m) = a.[s(m)]$ for each $a \in A,$ $s\in \Gamma(\mathcal{E})$ and $m\in M.$
The  $A$-product is given by
$$(s,t)_{\Gamma(\mathcal{E})} = \int_{ M} (s,t)_{\mathcal{A}} |\mbox{vol}_g|.$$
 
We often omit writing the reference to the bundle, and use the symbol $(,)_{\Gamma}$ instead of $(,)_{\Gamma(\mathcal{E})}.$ 
Let us notice that in the formula for the $A$-product $(,)_{\Gamma}$,  any absolutely convergent integral on Banach valued functions may be considered.
We take the so called Bochner integral for definiteness. 

{\bf Lemma 4:} If $(M,g)$ is a compact Riemannian manifold and $\mathcal{E}$ is an $A$-Hilbert bundle over $M,$ then
$\mathcal{E}'=TM \otimes \mathcal{E}$ and $\mathcal{E}''=T^*M \otimes \mathcal{E}$ are $A$-Hilbert bundles as well.

{\it Proof.} Let us set $a.(v \otimes c) = v \otimes a.c$ and $a.(\alpha \otimes c)=  \alpha \otimes a.c$ for any 
$a\in A,$ $c\in \mathcal{E}_m,$ 
$v\in T_mM,$  $\alpha \in T^*_mM$ and $m\in M.$ Further, set
$(u \otimes c, v \otimes d)_m = g_m(u,v)(c, d)_m \in A$ and $(\alpha \otimes c, \beta \otimes d)_m = g(\alpha, \beta)(c,d)_m$
for $c,d \in \mathcal{E}_m,$   $u, v \in T_mM$  and $\alpha, \beta \in T^*_mM,$ $m\in M.$
It is straightforward to verify that these relations define $A$-Hilbert bundles. $\Box$

In what follows, when given an $A$-Hilbert bundle $\mathcal{E},$ we always consider the bundles $\mathcal{E}'$ and $\mathcal{E}''$ with the $A$-Hilbert bundle 
structure described in Lemma 4.

{\bf Definition 1:} Let $p: \mathcal{E}\to M$ be an $A$-Hilbert bundle. We call a map $\nabla: \Gamma(\mathcal{E}) \to \Gamma(T^*M \otimes \mathcal{E})$
covariant derivative if for each function $f\in \mathcal{C}^{\infty}(M)$ and sections $s_1,s_2 \in \Gamma(\mathcal{E}),$ we have
\begin{eqnarray*}
\nabla (s_1 + s_2) &=& \nabla s_1 + \nabla s_2 \\
\nabla (f s_2) &=& df \otimes s_2 + f\nabla s_2
\end{eqnarray*}
For each $X \in \mathfrak{X}(M)=\Gamma(TM)$ and $s \in \Gamma(\mathcal{E}),$ we denote by $\nabla_X s$ the insertion $\iota_X(\nabla s) = (\nabla s)(X)$ of
$X$ into $\nabla s.$


Any covariant derivative $\nabla$ in an $A$-Hilbert bundle $\mathcal{E}$ induces the exterior covariant derivatives
$d_k^{\nabla}: \Gamma(\bigwedge^{k}T^*M \otimes \mathcal{E}) \to 
\Gamma(\bigwedge^{k + 1}T^*M \otimes \mathcal{E})$ by the formula
$$d^{\nabla}_k(\alpha \otimes s)  = d\alpha \otimes s + (-1)^{k}\alpha \wedge \nabla s$$
where $\alpha \otimes s \in \Gamma(\bigwedge^{k}T^*M \otimes \mathcal{E})$
  and $k=0, \ldots, \mbox{dim } M.$ To non-homogeneous elements, we extend the exterior covariant derivative by  linearity.
Further, we set $d_X^{\nabla} = \iota_X d^{\nabla}$ for any $X\in \mathfrak{X}(M).$

Let $\mathcal{E}, \mathcal{F}$ be $A$-Hilbert bundles over $M.$ Suppose that $\mathfrak{d}:\Gamma(\mathcal{E}) \to \Gamma(\mathcal{F})$ is a differential operator.
Then it is known that the symbol $\sigma: \mathcal{E}'' \to \mathcal{F}$ of $\mathfrak{d}$ is an adjointable $A$-Hilbert bundle morphism and 
that $\mathfrak{d}$ itself is adjointable as a pre-Hilbert $A$-module homomorphism.
Notice that we consider finite order operators only. 
Further, for each $t \in \mathbb{Z},$ one may define an $A$-product $(, )_t$ on $\Gamma(\mathcal{E}).$ The $A$-modules $\Gamma(\mathcal{E})$ equipped with
$(,)_t$ are pre-Hilbert $A$-modules.  Let us denote the norm associated to $(,)_t$  by $| \mbox{ } |_t.$ Notice that $(,)_0=(,)_{\Gamma(\mathcal{E})}.$
 The spaces 
$W^t(\mathcal{E})$   are defined as completions of $\Gamma(\mathcal{E})$ with respect to the norms $| \mbox{ } |_t.$  
Because the shape of the formulas for  $(,)_t$  is the same as the one for the classical Sobolev spaces, one may call the Hilbert $A$-module
$W^t(\mathcal{E})$ the Sobolev completion (of the pre-Hilbert $A$-module $(\Gamma(E), (,)_t)$). Let us notice that for any $t\in \mathbb{Z},$ each $A$-differential operator
$\mathfrak{d}:\Gamma(\mathcal{E}) \to \Gamma(\mathcal{F})$
 has a continuous extension to $W^t(\mathcal{E})$ and that this extension is unique. For these results, see Solovyov, Troitsky \cite{FM}.

We use the symbols $\xi^g$ and $\mbox{ext}_i^{\xi}$  introduced in section 2 to denote the $g$-dual vector field and the exterior multiplication by a differential $1$-form
$\xi$  also in the case of a Riemannian manifold $(M,g).$ 

{\bf Theorem 5:} Let $(M,g)$ be a compact Riemannian manifold, $p: \mathcal{E} \to M$ be an $A$-Hilbert bundle and $\nabla$ be a 
covariant derivative in $\mathcal{E}.$ Then $d^{\nabla}_i$ is an $A$-differential operator of order at most one. For each $\xi \in \Gamma(T^*M),$ 
the symbol $\sigma_i$ of $d_i^{\nabla}$ is given
by $\sigma_i^{\xi} = \mbox{ext}_i^{\xi} \otimes \mbox{Id}_{\mathcal{E}}$ and its adjoint satisfies $(\sigma_i^{\xi})^*  = \iota_{\xi^g} \otimes \mbox{Id}_{\mathcal{E}}.$

{\it Proof.}   For any function $f$ on $M$ and any section 
$\psi$ of $\bigwedge^i T^*M \otimes \mathcal{E},$ we get $d^{\nabla}_i(f\psi) - fd_i^{\nabla}\psi = (df) \wedge \psi +fd_i^{\nabla}\psi - f d_i^{\nabla}\psi=
df \wedge \psi$
which shows that the exterior covariant derivative $d_i^{\nabla}$ is a differential operator at most of first order.
For $\xi \in \Gamma(T^*M),$ let us compute the symbol  $\sigma_i^{\xi}$
of $d_i^{\nabla}.$ It is sufficient to work locally. Using the previous formula for $df =\xi$ and $\psi = \alpha \otimes s
\in \Gamma(\bigwedge^i T^*M \otimes \mathcal{E}),$ we get
$\sigma^{\xi}_i(\alpha \otimes s) = \xi \wedge \alpha \otimes s.$  In particular, the symbol acts on the form part only.
Because the adjoint of the wedge multiplication by a differential form $\xi$ is the interior product by the dual vector field $\xi^g,$ 
we get the formula
$(\sigma_i^{\xi})^* (\alpha \otimes s) = \iota_{\xi^g}\alpha \otimes s.$
$\Box$

{\bf Remark:} From the proof of the previous theorem, we see that the symbol of $d_i^{\nabla}$ is an adjointable 
homomorphism between the $A$-Hilbert bundles $(\bigwedge^i T^*M \otimes \mathcal{E})\otimes T^*M$ and $\bigwedge^{i+1}T^*M \otimes \mathcal{E}.$

\subsection{De Rham complex twisted by the oscillatory representation}

Let $A$ be a unital  $C^*$-algebra and  $\mathcal{E} \to M$ be an  $A$-Hilbert bundle over a compact manifold $M.$
 See Fomenko, Mishchenko \cite{FM}.

Let $(p_i:\mathcal{E}^i \to M)_{i \in \mathbb{N}_0}$ be a sequence of $A$-Hilbert bundles and
$\mathfrak{d}^{\bullet}=(\mathfrak{d}_i, \Gamma(\mathcal{E}^i))_{i\in \mathbb{N}_0}$ be a complex of $A$-differential operators.
We call such a  complex {\it $A$-elliptic} if out of the zero section of the cotangent bundle, the symbol sequence $\sigma^{\bullet}$ of $\mathfrak{d}^{\bullet}$
forms an exact complex in the category of $A$-Hilbert bundles. 
Notice that if the $A$-Hilbert bundles are vector bundles associated to a principal bundle,
 it is sufficient to demand the exactness of the symbol sequence on the level of the fibers, i.e., 
in the category of Hilbert $A$-modules.

Suppose  that for each $i\in \mathbb{N}_0,$ $\mathcal{E}_i \to M$ is a finitely generated projective $A$-Hilbert bundle over a compact manifold $M$, i.e., that for each $i,$ the fiber of 
$\mathcal{E}_i$ is  such  a Hilbert $A$-module. (See Solovyov, Troitsky \cite{ST} for a definition of a finitely generated projective Hilbert $A$-module.)
To any complex $\mathfrak{d}^{\bullet}=(\mathfrak{d}_i, \Gamma(\mathcal{E}^i))_{i \in \mathbb{N}_0}$ of differential operators, one may
 consider the sequence of its associated Laplacians
$\triangle_i = \mathfrak{d}_{i-1}\mathfrak{d}^{*}_{i-1} + \mathfrak{d}^{*}_{i}\mathfrak{d}_i,$ $i \in \mathbb{N}_0.$ Let us denote the order of $\triangle_i$
by $r_i.$ In Kr\'ysl \cite{Krysl} (Theorem 11), the following implication is proved.
If for each $i \in \mathbb{N}_0$ the continuous extension 
of the Laplacian $\triangle_i$ to $W^{r_i}(\mathcal{E}^i)$ has closed image, then the cohomology groups
$$H^i(\mathfrak{d}^{\bullet},A) = \frac{\mbox{Ker}(\mathfrak{d}_i:\Gamma(\mathcal{E}^i) \to \Gamma(\mathcal{E}^{i+1}))}{\mbox{Im}(\mathfrak{d}_{i-1}:\Gamma(\mathcal{E}^{i-1}) 
\to \Gamma(\mathcal{E}^{i}))}$$
of $\mathfrak{d}^{\bullet}$ are finitely generated projective $A$-modules and Banach spaces as well. 
The norm considered on the cohomology groups $H^i(\mathfrak{d}^{\bullet},A)$ is the quotient norm derived from the (possibly non-complete) norm 
$|\mbox{ }|_{0}$  on the smooth sections $\Gamma(\mathcal{E}^i).$

Now, let us focus to the specific case of symplectic manifolds and the de Rham complex tensored by the basic oscillatory module.
Let $(M^{2n},\omega)$ be a symplectic manifold.
In a similar way as in Riemannian geometry, one may introduce the notion of a "spin" structure over $(M,\omega),$ the so called {\it metaplectic structure}. 
See Habermann, Habermann \cite{KH} for a definition.
Suppose that $(M,\omega)$ possesses a metaplectic structure and denote it by $\tilde{\mathcal{G}}.$ In particular, 
$\tilde{\mathcal{G}}$ is a principal $\tilde{G}$-bundle over $M,$ where $\tilde{G}$ denotes  the metaplectic group.
Let $\mathcal{C}^{\bullet}$ denote the vector bundle associated to the principal bundle $\tilde{\mathcal{G}}$ via the representation
$\rho: \tilde{G} \to \mbox{Aut}(C^{\bullet}),$ i.e., $\mathcal{C}^{\bullet} = \tilde{\mathcal{G}}\times_{\rho} C^{\bullet}.$
Especially, $\mathcal{C}^0 = \tilde{\mathcal{G}}\times_{\rho} H$ is the so called {\it basic oscillatory bundle} which we denote by 
$\mathcal{H}$ here. In Habermann, Habermann \cite{KH} this bundle is called the symplectic spinor bunde. See also Kostant \cite{Kostant}.
Let $\nabla$ be a symplectic connection on $(M,\omega),$ i.e., $\nabla$ is a covariant derivative in $TM \to M$ preserving the symplectic form $\omega.$
We  allow the connection to have a non-zero torsion.
Let us denote a lift of this connection to $\tilde{\mathcal{G}}$ by $\tilde{\omega}^{\nabla}.$
Lifting $\tilde{\omega}^{\nabla}$ to $\mathcal{C}^{0},$ we get a covariant derivative $\nabla^H$ in the symplectic spinor bundle.
This covariant  derivative gives rise to a sequence $d_H^{\bullet}=(d^{\nabla^H}_i, \Gamma(\mathcal{C}^i))_{i=-1}^{2n+1}.$ 
Because $C^{i}$ is a Hilbert $A$-module (Lemma 3), the bundle $\mathcal{C}^i$ is an $A$-Hilbert bundle, where $A=\mbox{End}(H).$

{\bf Theorem 6:} Let  $(M^{2n},\omega)$  be a compact symplectic manifold which admits a  metaplectic structure, and $\nabla$ be a flat symplectic connection.
If the continuous extension to the Sobolev completions $W^2(\mathcal{C}^i)$ of each of the associated Laplacians
$\triangle_i$  has closed images, then the cohomology groups $H^i(d^{\bullet}_H,A)$ are Banach vector spaces and finitely generated projective $A$-modules.

{\it Proof.} Due to Lemma 3, the bundle $\mathcal{C}^{\bullet} \to M$ is a finitely generated projective $A$-Hilbert bundle.
Due to  Theorem 5, the symbol $\sigma_i$ of $d_i^{\nabla^H}$ is given by $\sigma_i^{\xi}  = \mbox{ext}_i^{\xi} \otimes \mbox{Id}_{\mathcal{H}}.$
 Thus the exactness of $(\sigma^{\xi}_i)_{i=-1}^{2n+1}$ is equivalent to the exactness of $\mbox{ext}^{\bullet}_{\xi}.$ The Cartan lemma (Lemma 4) implies that
 $(\sigma_{i}^{\xi}, \mathcal{C}^i)_{i=-1}^{2n+1}$
is exact and thus, $d^{\bullet}_H$ is an $A$-elliptic complex. Therefore Theorem 11 in \cite{Krysl} (mentioned above) may be applied and the conclusions on the cohomology groups
 $H^i(d^{\bullet}_H,A)$ follow.
$\Box$

The assumption on the  images of the extensions of the Laplacians seems to be unpleasant and we would like to conjecture
that it is satisfied in the case of the de Rham complex twisted by the basic oscillatory module.

\end{document}